\documentclass[sn-mathphys,iicol]{sn-jnl}

\usepackage[capitalise]{cleveref}
\usepackage[theme=default-plain]{jlcode}

\jyear{2023}%
\begin{document}

\title[Automating Adjoints]{Automating Steady and Unsteady Adjoints: Efficiently Utilizing Implicit and Algorithmic Differentiation}

\author*[1,2]{\fnm{Andrew} \sur{Ning}}\email{aning@byu.edu}

\author[1]{\fnm{Taylor} \sur{McDonnell}}

\affil[1]{\orgdiv{Mechanical Engineering}, \orgname{Brigham Young University}, \orgaddress{\street{360 EB}, \city{Provo}, \postcode{84602}, \state{UT}, \country{USA}}}

\affil[2]{\orgdiv{Wind Energy}, \orgname{National Renewable Energy Laboratory}, \orgaddress{\street{15013 Denver West Pkwy}, \city{Golden}, \postcode{80401}, \state{CO}, \country{USA}}}

\abstract{Algorithmic differentiation (AD) has become increasingly capable and straightforward to use.  However, AD is inefficient when applied directly to solvers, a feature of most engineering analyses. We can leverage implicit differentiation to define a general AD rule, making adjoints automatic.  Furthermore, we can leverage the structure of differential equations to automate unsteady adjoints in a memory efficient way.  We also derive a technique to speed up explicit differential equation solvers, which have no iterative solver to exploit.  All of these techniques are demonstrated on problems of various sizes, showing order of magnitude speed-ups with minimal code changes.  Thus, we can enable users to easily compute accurate derivatives across complex analyses with internal solvers, or in other words,  automate adjoints using a combination of AD and implicit differentiation.}


\maketitle

\section{Introduction}

Implicit equations require the use of a solver and are a feature of most engineering analyses.  The iterative nature of solvers not only complicates the analysis, but also makes it more difficult to obtain accurate derivatives efficiently.  This is especially true for solvers with long run times.  Implicit analytic methods \cite{Martins2022}, commonly referred to as direct and adjoint methods, help address this problem and have been in use for a long time, starting in the field of optimal control \cite{Bryson1996} and later applied to structural \cite{Arora1979} and aerodynamic disciplines \cite{Pironneau1974,Jameson1988}, as well as coupled systems \cite{Martins2005}. In this paper we will use the term adjoint for brevity, but really mean either a direct (forward) or adjoint (reverse) formulation, as appropriate for the given problem.

An adjoint uses implicit differentiation to analytically obtain the desired total derivatives from partial derivatives.  These partial derivatives are much easier to obtain as they do not depend on solver iterations.  The accuracy of the method then hinges on the accuracy of the partial derivatives.  When exact methods are used to obtain the partials, the total derivatives can be exact.

While adjoint methods are well known and used, implementing one for a given problem can be incredibly time consuming.
AD, on the other hand, provides exact partials automatically, and avoids the need to specify a connectivity graph, but is not efficient across solvers.  The iterations of the solver, can be unrolled, and AD applied to this unrolled set of operations, but this is usually less efficient and less accurate. In many cases it is not even possible as the internals of the particular solver may not be overloaded for AD.  This process of applying AD to the unrolled set is sometimes called \emph{direct AD}.  Each implementation of AD is also limited to a specific programming or modeling language.

AD is based on \emph{rules}, which define the exact derivatives of a given numerical operation (e.g., the derivative of multiplication, or the cosine function).  AD rules, in combination with the chain rule, enable automated differentiation of arbitrary codes \cite{Martins2022}.  If we provide rules for solvers, then we can enable usage that is as simple as direct AD, but is overloaded to make use of adjoints automatically.  Because AD is no longer used within the solvers, we do not need AD compatibility within solvers, in fact we could utilize solvers from other programming languages.  Additionally, we find that the same methodology is easily extended to allow for custom Jacobians (or vector-Jacobian products).  In other words, the AD interface can now easily incorporate other programming languages or differentiation strategies (including mixed-mode AD).

Using implicit differentiation to create an AD rule for solvers (also known as a super node \cite{Margossian2019}) is well known and has existed for some time \cite{Schachtner1993,Giering1998,Christianson1998,Bartholomew-Biggs2000,Bell2008,Blondel21}.  However, its implementation has become much easier in recent years as AD has become much more widely available and user-friendly.  This is especially true in the Julia programming language \cite{Bezanson2017} where efficient operator overloaded AD has been a core part of the language since its inception.  While many languages provide AD support to varying degrees, operator-overloaded AD in Julia is essentially as straightforward as finite differencing and as efficient as source-code transformation.

While the adjoint can be applied to any set of implicit equations, the implementation associated with time-based problems is often given the name \emph{unsteady adjoint}. An unsteady adjoint is applicable to differential equations in which a given state only depends on prior states (and usually only a small subset of prior states).  Time is the most common example, but spatial or other derivatives also sometimes follow this form.  Ordinary differential equations (ODEs), many differential algebraic equations, and some solution methods for partial differential equations fall under this form.  For these problems, especially time based ones, computational and memory requirements can be extensive requiring special consideration. As ODE problems grow larger, a monolithic application of AD is prohibitively slow.

Various approaches exist to address this problem, such as solving a larger ODE system with the sensitivities added as additional state variables \cite{Carpenter2015}, using continuous adjoint approaches (discussed in the next paragraph) \cite{Chen2018,Ma2021}, using various levels of checkpointing on discrete approaches \cite{Zhuang2020,Zhang2022}, or using the unified derivative equation approach \cite{Hwang2018} where the partials and graph structure must be provided and then total derivatives can be computed from a linear solve \cite{Falck2021}.

Like regular adjoints, unsteady adjoints come in both continuous and discrete varieties. In the former, implicit differentiation is applied to the continuous, or a semi-discrete, set of differential equations, and discretization occurs at the end \cite{Bradley2010}.  Whereas in the latter approach, we discretize the differential equations first, then apply implicit differentiation \cite{Betancourt2020}.  Our focus is on the discrete form, which is not always as efficient \cite{Ma2021}, but is consistent with the implementation of the differential equations (the continuous approach is only consistent in the limit as the discretization is refined).  The inconsistency of the continuous approach can sometimes mislead the optimizer \cite{Peter2010}, whereas the  derivatives from the discrete approach are generally more accurate and stable \cite{Onken2020}.

In this first section we address efficient derivative computation for solvers of nonlinear functions.  Any implicit function that we wish to solve can be classified under the nonlinear case.  However, for some functions we can obtain derivatives more efficiently than just treating it as a general nonlinear problem.  For example, analytic solutions exist for linear systems in forward mode and reverse mode \cite{Giles2008}.  Additionally, eigenvalue and eigenvector problems also have analytic solutions in forward mode \cite{Wilkinson1965,Magnus1998,Leeuw2007} and reverse mode \cite{He2022} (the latter with some novel solutions for eigenvectors).  

Next, we show how to obtain derivatives efficiently for differential equations in forward and reverse mode using implicit solvers.  In particular, we highlight how to compute discrete unsteady adjoints in a memory efficient manner, even with long time sequences, for any ODE solver.  We also demonstrate an approach for explicit ODEs that reduces the cost of generating a long computational graph.  Finally, we detail examples for these various categories of solvers, with variable size problems so that we can compare the scalability of the different methods.

We have implemented all the methodology shown here, including analytic overloads for linear solvers and eigenvalue solvers, in the open source package ImplicitAD.\footnote{\url{https://github.com/byuflowlab/ImplicitAD.jl}}  The math and implementation is fairly straightforward and lends itself well to simplified implementations for instructional purposes.

\section{Theory and Implementation}

Primal is used to refer to the original engineering output, and we will use partials to refer to the corresponding derivative(s).  In forward mode, this is often implemented in what is called a dual number, an object containing both primal and partial numbers (somewhat analogous to a complex number).

\subsection{Nonlinear Functions}
\label{sec:nonlinear}

We can write any solver, in residual form as:
\begin{equation}
    r(x, y(x)) = 0
    \label{eq:genresid}
\end{equation}
where $r$ is the residual, $y$ are the corresponding states (implicit outputs), and $x$ are the inputs.  For the purposes of AD we are after the Jacobian $dy/dx$ for this operation (the application of a solver to the above residuals).
Using the chain rule we find that:
\begin{equation}
\frac{dr}{dx} = \frac{\partial r}{\partial x} + \frac{\partial r}{\partial y} \frac{dy}{dx} = 0,
\label{eq:ddx}
\end{equation}
which we solve for the desired Jacobian:
\begin{equation}
    \frac{dy}{dx} = -\left( \frac{\partial r}{\partial y} \right)^{-1}\frac{\partial r}{\partial x}
    \label{eq:jacobian1}
\end{equation}
Note that the solution depends only on partial derivatives at the solution of the residual equations, and thus not on the solver path, which would be traversed in direct AD.
In the typical adjoint derivation, an output function of interest is also defined, which is an explicit function of $x$ and $y$. However, since we are using AD those derivatives are propagated automatically (either explicitly with regular AD or implicitly using the methods discussed here).  Thus, we need only the intermediate step of computing $dy/dx$ for this implicit function, i.e., defining an AD rule for a generic nonlinear solver.

\subsubsection{Forward Mode Implicit AD}

If we are using forward mode algorithmic differentiation, once we reach this implicit function we already know the derivatives of input $x$ with respect to our overall input variables of interest ($\eta$).  Note that the inputs $\eta$ apply to the entire computational program, whereas $x$ are inputs to just this solver operation.  Using common AD nomenclature we call this derivative $\dot{x}$, shown below for one variable of interest:
\begin{equation}
\dot{x}_i \equiv \frac{\partial x_i}{\partial \eta}
\end{equation}
\Cref{eq:jacobian1} describes how to compute the Jacobian, but what we really want is the Jacobian vector product (JVP) that yields $\dot{y}$:
\begin{equation}
    \dot{y} = \frac{dy}{dx} \dot{x}
\end{equation}

We multiply both sides of our governing equation \cref{eq:jacobian1} by $\dot{x}$, and rearrange:
\begin{equation}
\frac{\partial r}{\partial y} \dot{y} = -\frac{\partial r}{\partial x} \dot{x}
\end{equation}
We focus first on the right-hand side, which is a JVP.
\begin{equation}
b = -\frac{\partial r}{\partial x} \dot{x}
\label{eq:jvp}
\end{equation}
We can compute the result without actually forming the matrix $\partial r / \partial x$, in an efficient manner.  The JVP is a weighted sum of the columns
\begin{equation}
b = \sum_i \frac{\partial r}{\partial x_i} (-\dot{x}_i)
\end{equation}
where $-\dot{x}_i$ are the weights.  Normally we pick off one of these partials in forward mode AD (i.e., a column of the Jacobian) one at a time with different seeds $e_i = [0, 0, 1, 0]$, but if we set the initial seed to $-\dot{x}_i$ then the partial derivatives we compute are precisely this weighted sum.  Thus, we can compute the desired JVP in one forward-mode pass independent of the size of $x$ or $y$.

In our implementation we compute the JVP by simply evaluating the residual with $y$ as a constant input and $\dot{x}$ as a dual number input.  We then extract the resulting derivative portion of the resulting dual numbers: \jlinl{b = -partials(r(xdual, y))}.

We now just need to compute the square Jacobian $\partial r / \partial y$ and solve the linear system.
\begin{equation}
\frac{\partial r}{\partial y} \dot{y} = b
\label{eq:forwardlin}
\end{equation}
Often this Jacobian is already available from the solver $r(x, y(x))$.  If not, because $\partial r / \partial y$ is square, forward mode AD is usually preferable for computing these partial derivatives. If we know the structure of the Jacobian, we would want to provide an appropriate factorization.  Sometimes, this Jacobian is sparse, and so using graph coloring is desirable.  Also, for large systems we can use matrix-free methods (e.g., Krylov methods) to solve the linear system without actually constructing the matrix.

In this forward mode we see that we must solve a linear system for each input $\eta$, thus the cost scales with the size of our overall inputs.

\subsubsection{Reverse Mode Implicit AD}
\label{sec:revnonline}
If using reverse mode AD, once we reach this implicit function we have the derivatives of our overall outputs of interest ($\xi$) with respect to our intermediate outputs $y$.  Using the common AD nomenclature (shown for one variable of interest):
\begin{equation}
\bar{y}_i \equiv \frac{\partial \xi}{\partial y_i}
\end{equation}
We need to propagate backwards to get $\bar{x}$.

Again, we are not actually interested in the Jacobian shown in \cref{eq:jacobian1}, but rather in the vector Jacobian product (VJP) that yields $\bar{x}$:
\begin{equation}
    \bar{x}^T = \bar{y}^T \frac{dy}{dx}
\end{equation}
Multiplying both sides of our governing equation (\cref{eq:jacobian1}) by $\bar{y}^T$ gives:
\begin{equation}
\bar{x}^T = - \bar{y}^T \left(\frac{\partial r}{\partial y} \right)^{-1} \left(\frac{\partial r}{\partial x} \right)
\end{equation}
We call the first two terms on the right-hand side  $\lambda^T$ (without the minus sign), and thus need to solve the linear equation:
\begin{equation}
\left(\frac{\partial r}{\partial y} \right)^T \lambda = \bar{y}
\label{eq:revlin}
\end{equation}
Like the previous section, we need to provide or compute the square Jacobian $\partial r / \partial y$.  Again, if the Jacobian (or an appropriate factorization or matrix-free operator) is not already available from the outer solver, then forward mode AD is usually preferable for these partials since the Jacobian is square.

Once we solve for $\lambda$ we see that we get our desired derivatives from the VJP:
\begin{equation}
    \bar{x}^T = - \lambda^T \frac{\partial r}{\partial x}
    \label{eq:vjp}
\end{equation}
A VJP is desirable as we can compute it very efficiently without actually forming the matrix $\partial r / \partial x$.  Note that the VJP is a weighted sum of the Jacobian's rows
\begin{equation}
\bar{x} = \sum_i -\lambda_i \frac{\partial r_i}{\partial x}
\end{equation}
Normally we pick off these partials in reverse mode AD one row at a time, but if we set the initial seed (called the adjoint vector in AD) to $-\lambda_i$ then the partial derivatives we compute are precisely this weighted sum.  Thus, we can compute the desired VJP in one pass independent of the size of $x$ or $y$.

In our implementation we solve \cref{eq:vjp} by computing the gradient of $h(x) = -\lambda^T r(y, x)$ with reverse-mode AD.  Since $\lambda$ is a constant, this gives us the desired VJP, and is efficiently obtained with reverse-mode AD since the above dot product produces a single scalar output.

Again, the main cost is in solving a linear system and so the overall cost scales with the size of $\xi$, the number of overall outputs.

\subsubsection{Fixed Point Problems}

Some nonlinear problems can be written in fixed point form:
\begin{equation}
y = f(x, y(x))
\end{equation}
We can rewrite this back in the general residual form (\cref{eq:genresid}) as:
\begin{equation}
    r(x, y) = f(x, y) - y = 0
\end{equation}
We can then reuse everything from the previous section.  However, we can avoid creating the trivial residual and directly obtain the partials from $f$ as:
\begin{equation}
\frac{\partial r}{\partial x} = \frac{\partial f}{\partial x}
\end{equation}
and
\begin{equation}
\frac{\partial r}{\partial y} = \frac{\partial f}{\partial y} - I
\end{equation}
where $I$ is the identity matrix.

\subsection{Differential Equations}

A set of ordinary differential equations can be expressed as a set of residual expressions ($r$), starting with an initialization, followed by a sequence of discrete integration steps from time $t_0$ to $t_{n}$ that update the states $y$ (with $x$ as a set of fixed inputs).
\begin{equation}
r(x, y, t) = \begin{bmatrix}
r_0(x, y_0, t_0) \\
r_1(x, y_{0 \cdots 1}, t_{0 \cdots 1}) \\
r_2(x, y_{0 \cdots 2}, t_{0 \cdots 2}) \\
r_3(x, y_{0 \cdots 3}, t_{0 \cdots 3}) \\
\vdots \\
r_n(x, y_{0 \cdots n}, t_{0 \cdots n}) \\
\end{bmatrix} = 0
\label{eq:oderesids}
\end{equation}
While the residuals do have time dependence, time is not typically a function of $x$ and so our original chain-rule based equation (\cref{eq:ddx}) is unchanged (if time is dependent on $x$ we can use the same equations by treating time as an augmented state variable).
The process in \cref{sec:nonlinear} is then applied to the residuals.

\subsubsection{Forward Mode}

Recall that for forward mode, the derivatives of the outputs $\dot{y}$ may be found by first computing the Jacobian vector product shown in \cref{eq:jvp}, then solving the linear system shown in \cref{eq:forwardlin}.

Expanded for our residuals, this system of equations takes the form:
\begin{equation}
\begin{bmatrix}
{\frac{\partial r_0}{\partial y_0}} & & & & \\
 \frac{\partial r_1}{\partial y_0} & \frac{\partial r_1}{\partial y_1} & & & \\
 \vdots & \vdots & \ddots & & \\
  \frac{\partial r_{n-1}}{\partial y_0} & \frac{\partial r_{n-1}}{\partial y_1} & \cdots & \frac{\partial r_{n-1}}{\partial y_{n-1}} & \\
 \frac{\partial r_{n}}{\partial y_0} & \frac{\partial r_n}{\partial y_1} & \cdots & \frac{\partial r_n}{\partial y_{n-1}} & \frac{\partial r_n}{\partial y_n}
\end{bmatrix}
\begin{bmatrix}
  \dot{y}_0 \\
  \dot{y}_1 \\
  \vdots & \\
  \dot{y}_{n-1} \\
  \dot{y}_n
\end{bmatrix}
 =
\begin{bmatrix}
  b_0 \\
  b_1 \\
  \vdots \\
  b_{n-1} \\
  b_n
\end{bmatrix}
\label{eq:fbigmatrix}
\end{equation}
where
\begin{equation}
    b_i = -\frac{\partial r_i}{\partial x} \dot{x}
\end{equation}

We can solve this problem iteratively starting from the first time step:
\begin{equation}
\frac{\partial r_0}{\partial y_0} \dot{y}_0 = -\frac{\partial r_0}{\partial x} \dot{x}
\label{eq:f_ic}
\end{equation}
For the $i^{\text{th}}$ set of equations we have:
\begin{equation}
    \sum_{k=0}^{i} \frac{\partial r_i}{\partial y_k} \dot{y}_k = b_i
\end{equation}
Since all the derivatives with respect to $i-1$ and lower are already known from previous lines we move those to the right hand side.
\begin{equation}
    \frac{\partial r_i}{\partial y_i} \dot{y}_i = b_i - \sum_{k=0}^{i-1} \frac{\partial r_i}{\partial y_k} \dot{y}_k
\end{equation}

In practice a given residual doesn't depend on all previous states, but just a small subset of prior time steps.  Thus, the above matrix is block banded (or block diagonal for one-step methods).  Below we'll assume dependence on the prior two states, if more states needed, the additional terms have the same form.  One-step methods are common and in those cases, we would just drop that last term.  After solving for $\dot{y}_0$ in \cref{eq:f_ic} we move forward in time to solve $\dot{y}_i$ from $i = 1 \ldots n$.
\begin{equation}
\frac{\partial r_i}{\partial y_i} \dot{y}_i = - \left[ \frac{\partial r_i}{\partial x} \dot{x} + \frac{\partial r_i}{\partial y_{i-1}} \dot{y}_{i-1} + \frac{\partial r_i}{\partial y_{i-2}} \dot{y}_{i-2}  \right]
    \label{eq:ode_f_iter}
\end{equation}
The entire right hand-side can be computed efficiently with a single Jacobian vector product.

Note that because of our reliance on AD, the above derivation is completely general, and thus applies to any set of differential equations that follows the structure in \cref{eq:oderesids}.

\subsubsection{Reverse Mode}

For reverse mode we first solve the linear system \cref{eq:revlin}, which we expand out for our system of equations:
\begin{equation}
\begin{bmatrix}
{\frac{\partial r_0}{\partial y_0}}^T & \frac{\partial r_1}{\partial y_0}^T & \cdots & \frac{\partial r_{n-1}}{\partial y_0}^T & \frac{\partial r_{n}}{\partial y_0}^T \\
 & \frac{\partial r_1}{\partial y_1}^T & \cdots & \frac{\partial r_{n-1}}{\partial y_1}^T & \frac{\partial r_n}{\partial y_1}^T \\
 &  & \ddots & \vdots & \vdots \\
  & & &  \frac{\partial r_{n-1}}{\partial y_{n-1}}^T & \frac{\partial r_n}{\partial y_{n-1}}^T\\
 & & & & \frac{\partial r_n}{\partial y_n}^T
\end{bmatrix}
\begin{bmatrix}
  \lambda_0 \\
  \lambda_1 \\
  \vdots & \\
  \lambda_{n-1} \\
  \lambda_n
\end{bmatrix}
 =
\begin{bmatrix}
  {\bar{y}_0} \\
  {\bar{y}_1} \\
  \vdots \\
  {\bar{y}_{n-1}} \\
  {\bar{y}_n}
\end{bmatrix}
\label{eq:rbigmatrix}
\end{equation}

This time we see we can start from the final time step to solve for $\lambda_n$:
\begin{equation}
    \frac{\partial r_n}{\partial y_n}^T \lambda_n = \bar{y}_n
    \label{eq:r_ic}
\end{equation}
and now work backwards in time:
\begin{equation}
    \frac{\partial r_i}{\partial y_i}^T \lambda_i = \bar{y}_i - \sum_{k=i+1}^n \frac{\partial r_k}{\partial y_i}^T \lambda_k
\end{equation}
Again, we only need to retain one or two terms in practice. So after solving for $\lambda_n$ in \cref{eq:r_ic} we work backwards in time solving $\lambda_i$ for $i = (n-1) \ldots 0$.
\begin{equation}
    \frac{\partial r_i}{\partial y_i}^T \lambda_i = \bar{y}_i - \frac{\partial r_{i+1}}{\partial y_i}^T \lambda_{i+1} - \frac{\partial r_{i+2}}{\partial y_i}^T \lambda_{i+2}
    \label{eq:ode_adj}
\end{equation}
Note that the last term is not included for the $i = n-1$ case.

After solving the linear system we compute the vector Jacobian product shown in \cref{eq:vjp}.
To reduce memory requirements, we can use the elements of $\lambda$ as they are computed to find their contribution to $\bar{x}$ using the following expression
\begin{equation}
\bar{x}^T = \sum_{i=0}^n -\lambda_i^T \frac{\partial r_i}{\partial x}
    \label{eq:xupdate}
\end{equation}
where we run the summation in reverse.

To further decrease computational expenses, the vector Jacobian products $\lambda_i^T \frac{\partial r_i}{\partial x}$, from the above expression, and  $\lambda_i^T \frac{\partial r_i}{\partial y_{i-1}}$ from the right hand side of \cref{eq:ode_adj}, may be computed simultaneously without forming the Jacobian using reverse-mode automatic differentiation.

In addition to an approach that is general, we also note that memory requirements are kept small.  This is particularly advantageous as memory often because a limiting problem for many discrete unsteady  adjoint implementations.  The matrix in \cref{eq:rbigmatrix} is never formed.  The Jacobians in the right hand-side of the \cref{eq:ode_adj,eq:xupdate} are not formed either (computed as vector Jacobian products as noted).  The only Jacobian we need to form is $\partial r_i / \partial y_i$, but this is a small matrix because it applies only to a single instance in time.  Furthermore, because \cref{eq:ode_adj} can be updated at each time instance we only need to save a single vector for the previous time step (which can then be replaced).  The only information that needs to be saved from the forward pass is the solution to the states ($y$), which one would need anyway even if derivatives were not of interest.  Thus, this approach is memory efficient even for differential equations that require a large number of time steps.

\subsubsection{Explicit Methods}
\label{sec:explicitode}

If an explicit method is applied to solving the differential equations then one can just use direct AD as applied to the entire ODE:
\begin{equation}
y_n = g(x, y_{0, \cdots n-1}, t_{0 \cdots n})
\end{equation}
This is fine for forward mode AD, but reverse-mode AD will increase in memory requirements rapidly.  Operator overloaded reverse mode AD requires not only storing a new data type like in forward mode, but we must also save the computational graph during the forward pass so that we know what sequence of operations to execute during the reverse pass \cite{Martins2022}.  This process is called taping, and we often refer to the saved state of the computational graph as the \emph{tape}. Thus, long computational sequences produces long tapes.

Instead, it may be beneficial to break up the problem into individual time steps, as before, except this time each line is explicit:
\begin{equation}
\begin{bmatrix}
y_0 = f_0(x, t_0) \\
y_1 = f_1(x, y_{0}, t_{0 \cdots 1}) \\
y_2 = f_2(x, y_{0 \cdots 1}, t_{0 \cdots 2}) \\
y_3 = f_3(x, y_{0 \cdots 2}, t_{0 \cdots 3}) \\
\vdots \\
y_n = f_n(x, y_{0 \cdots n-1}, t_{0 \cdots n}) \\
\end{bmatrix}
\end{equation}
We can rewrite this in the residual form of \cref{eq:oderesids} as:
\begin{equation}
    r(x, y, t) =
    \begin{bmatrix}
    f_0(x, t_0) - y_0\\
    f_1(x, y_{0}, t_{0 \cdots 1}) - y_1\\
    f_2(x, y_{0 \cdots 1}, t_{0 \cdots 2}) - y_2\\
    f_3(x, y_{0 \cdots 2}, t_{0 \cdots 3}) - y_3\\
    \vdots \\
    f_n(x, y_{0 \cdots n-1}, t_{0 \cdots n}) - y_n\\
    \end{bmatrix} = 0
    \end{equation}

Then, we can compute the residuals as:
\begin{equation}
\frac{\partial r_i}{\partial x} = \frac{\partial f_i}{\partial x}
\end{equation}
and
\begin{equation}
\frac{\partial r_i}{\partial y_j} = \begin{cases}
    -I & \text{ if } i = j \\
    \frac{\partial f_i}{\partial y_j} & \text{ if } i > j
\end{cases}
\end{equation}
where $I$ is the identity matrix.

In reverse mode, we substitute into \cref{eq:r_ic,eq:ode_adj} to obtain:
\begin{equation}
-\lambda_n = \bar{y}_n
\end{equation}
\begin{equation}
-\lambda_i = \bar{y}_i - \frac{\partial f_{i+1}}{\partial y_i}^T \lambda_{i+1} - \frac{\partial f_{i+2}}{\partial y_i}^T \lambda_{i+2}
\end{equation}
where again the last term is omitted for the $n-1$ case.  Thus, at time step $i$ (recall that we are running backwards) we compute the vector Jacobian products
\begin{equation}
\lambda_i^T \frac{\partial f_i}{\partial y_{i-1}}
\quad\text{ and }\quad
\lambda_i^T \frac{\partial f_i}{\partial y_{i-2}}
\end{equation}
and subtract them from the known values for $\lambda_{i-1}$ and $\lambda_{i-2}$ repectively.

We also see that each $\lambda_i$ is precisely the weights for the vector Jacobian product we need in \cref{eq:xupdate}.  At each time step, the three vector Jacobian products can be computed simultaneously as they use the same weights $\lambda_i$, and operate on the same function $f_i$.  Thus, as expected, we compute a sequence of vector Jacobian products starting from time $n$ and working backwards to time $0$.  However, this approach does not require recording a long tape as does a direct application of reverse-mode AD.  Because we only apply reverse-mode AD across each time step separately, and then analytically propagate to the previous time step, the tape is \emph{much} shorter.  This allows us to compute adjoints for long time-based simulations without incurring large memory penalties.  Additionally, because each time step calls the same function, just with different inputs, we can compile the tape and reuse it for each time step to significantly speed up computation time.  Note that some ODEs have branching behavior (e.g., if/else statements), and in those cases we cannot reuse the tape for each time step, and would instead prefer a source-to-source reverse-mode implementation.

\subsection{External Code}

While the focus of this work is on automating adjoints, the methodology easily allows for one to insert other custom derivatives within an AD chain.  A common use case would be with mixed programming languages.  In other words, by defining an AD rule for external function calls, AD can be applied in a direct manner.

We will refer to a generic explicit function, which would generally be in another languages as: $z = z(x)$.
For forward mode we use the chain rule, where $\eta$ is again our overall input variable of interests
\begin{equation}
\frac{d z_i}{d \eta_j} = \frac{d z_i}{d x_k} \frac{d x_k}{d \eta_j}
\end{equation}
or in our notation
\begin{equation}
\dot{z} = J \dot{x}
\label{eq:forwardexternal}
\end{equation}
where $\dot{x}$ is known at this stage of the AD chain and $\dot{z}$ is what we are solving for.  Ideally the external code can supply $J$, or Jacobian vector products, but if not, we could compute it with complex step or finite differencing if necessary.

We are actually after a Jacobian-Jacobian product, so if the Jacobian is not provided we have two options.  Either, we compute the Jacobian $J$ first (one column at a time), which scales with the dimension of $x$ (each index k from equation above).  Or, we can compute the Jacobian vector product for each $\eta$, which scales with the dimension of $\eta$ (over index j).  For example, with forward finite differencing we can compute the JVP by taking steps in the $\dot{x}$ direction (as opposed to steps in the coordinate directions) as:
\begin{equation}
\dot{z}_j \approx \frac{z(x + h \dot{x}_j) - z(x)}{h}
\end{equation}
We would compare the size of $x$ and $\eta$ and choose the smaller number of function calls.

In reverse mode we use the chain rule, where again $\xi$ are the overall output variables of interest:
\begin{equation}
\frac{d \xi_j}{d x_k} = \frac{d z_i}{d x_k} \frac{d \xi_j}{d z_i}
\end{equation}
or in our notation
\begin{equation}
\bar{x} = J^T \bar{z}
\label{eq:reverseexternal}
\end{equation}
where $\bar{z}$ is known at this stage of the AD chain (though generally given for just one value of $\xi$ at a time) and $\bar{x}$ is what we are computing.

In this case the user can again provide the Jacobian or the vector-Jacobian product (i.e., the gradient of $h(x) = \bar{z}^T z(x)$ ideally via reverse-mode AD), or we can fall back to complex step or finite differencing.  Finite differencing and complex step only work in a forward manner and so we cannot compute the VJP directly, the only option in reverse mode is to construct the Jacobian  and then multiply. If we had the whole Jacobian $\bar{z}$ available we could actually perform a JVP for each $z_i$, but the pullback operations in reverse AD provide the opposite.

\section{Examples}

AD is executed with the packages ForwardDiff.jl \cite{Revels2016} and ReverseDiff.jl\footnote{\url{https://github.com/JuliaDiff/ReverseDiff.jl}} for forward mode and reverse mode operator-overloaded AD respectively.  Various other reverse-mode AD options exist in Julia, some of which can offer unique performance advantages, depending on the problem.  ReverseDiff is chosen here for its robustness.

Performance timing is primarily used to show that implicit AD can offer significant speed advantages over direct AD, and to highlight general trends.  However, the relative advantages of the various methods will always be problem specific and one is encouraged to do their own benchmarking. Timings were performed on a laptop using an Apple M1 Pro chip with 16 GB RAM.

\subsection{Adjoint for a Nonlinear Solver}

We consider both a nonlinear problem that allows the number of states to be varied.  For both forward and reverse mode, all caches are preallocated (although this makes little difference in the timings), and the tape in reverse mode is precompiled (which improves speed by about a factor of 2 on this problem).

This example uses the $n$-dimensional Rosenbrock problem, but posed as a root-finding problem.  The residuals are:
\begin{equation}
    r(y) =
    \begin{cases}
-4 \alpha_1 y_1 (y_2 - y_1^2) - 2(1 - y_1) & i = 1\\
-4 \alpha_i y_i (y_{i+1} - y_i^2) \\
\quad - 2(1 - y_i)  \\
\quad + 2 \alpha_{i-1} (y_i - y_{i-1}^2) & i = 2 \ldots n-1\\
 2 \alpha_{n-1} (y_n - y_{n-1}^2) & i = n
\end{cases}
\label{eq:rosen}
\end{equation}
where $\alpha_i$ are parameters that are 100 in the original Rosenbrock problem.

We choose to use $n$ design variables evaluated at the point $x_i = 100$, where the coefficients are set as $\alpha_i = x_i$.  For outputs we just use the states $y$.  Thus, this problem has $n$ inputs, $n$ states, and $n$ outputs.  Because the number of inputs and outputs are equal, neither forward nor reverse mode offers a unique scaling advantage.

We solve problems of increasing size doubling $n$ from $n = 2, 4, 8, \ldots, 128$.  The nonlinear solver is the trust region method in NLsolve.jl \cite{nlsolve}.  For each size we compare forward mode and reverse mode using both direct AD and implicit AD, and also include a central finite differencing for comparison.  Even with central differencing the comparison is not quite fair as the finite differencing produces much less accurate derivatives, however it is still a useful reference point.

\begin{figure}[htbp]
\centering
\includegraphics[width=3.0in]{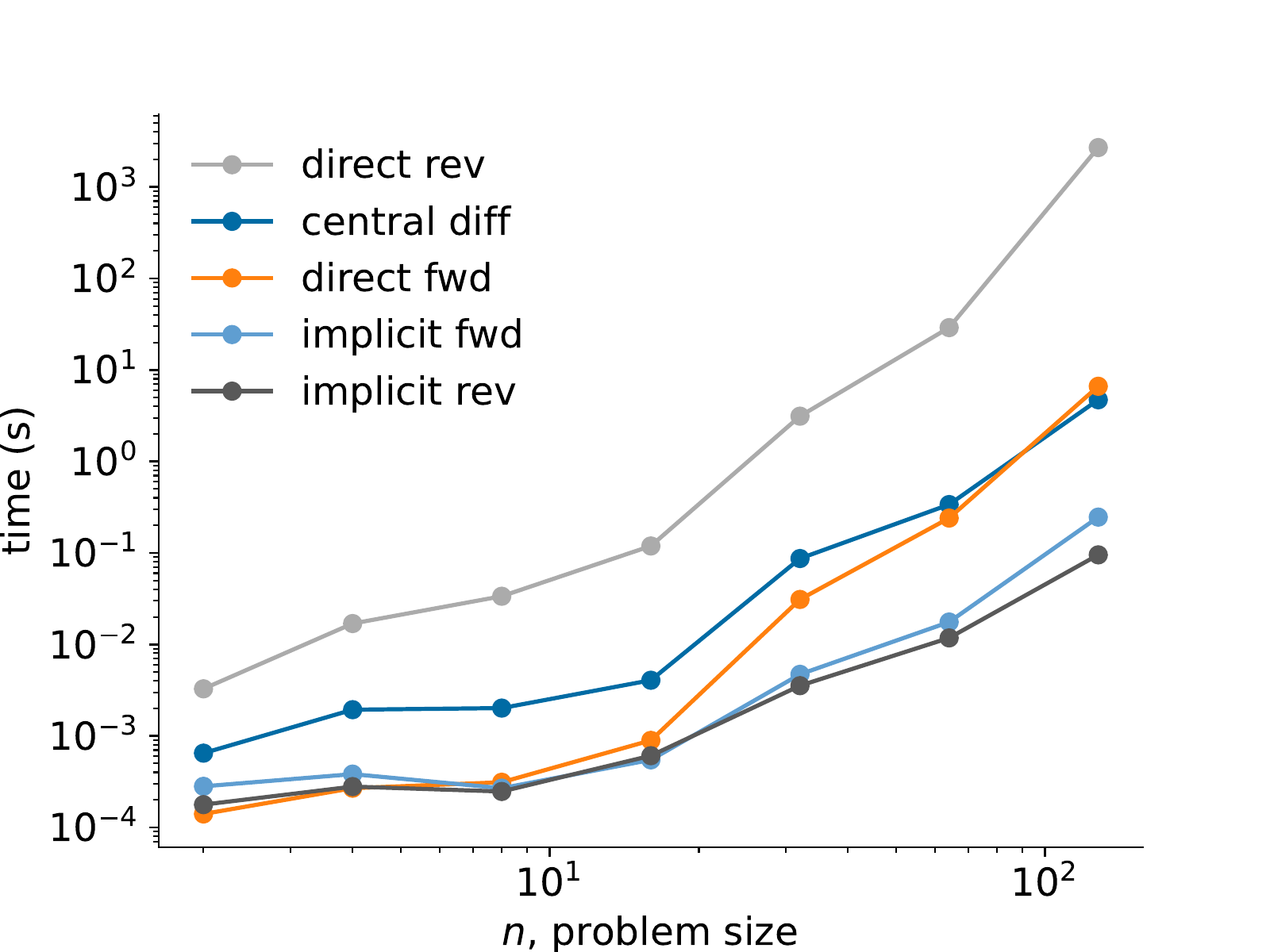}
\caption{A comparison in computing a dense Jacobian between direct forward direct, direct reverse, implicit forward, implicit reverse, and finite differencing for $n$ inputs and $n$ outputs.}
\label{fig:nonlinear-nn}
\end{figure}

In \cref{fig:nonlinear-nn} we show results for this case with $n$ inputs, $n$ states, and $n$ outputs.  The $y$ axis shows the median time to compute the Jacobian.  The median is computed across 10,000 samples (except for a handful of cases with relatively long compute times).  Because the values vary across orders of magnitude a log-log scale is utilized.

We notice that direct reverse AD is the slowest.  The problem dimensions didn't favor reverse mode to begin with, and creating the tape through the increasing number of solver iterations requires extra computations.  It is  particularly memory intensive since we are storing information for the reverse pass.

By contrast, we see that implicit reverse AD is the fastest (though not by much as compared to implicit forward AD).  In other words, just switching from direct to implicit takes us from the slowest method to the fastest. As discussed in the methodology, implicit differentiation bypasses the solver iterations---we no longer need to create a tape for the solver iterations.  Instead, we only need to apply reverse-mode AD to vector Jacobian products associated with just the residual function.  This only requires a very small tape.

With implicit differentiation, both forward and reverse mode are of comparable speeds, generally within a factor of 2 or less.  We expect this as the problem has equal numbers of inputs and outputs.  As we saw in the methodology, the main difference in the implicit forward versus implicit reverse, is the number of linear solves that are performed, which scales with the number of inputs or outputs respectively.  If that balance were shifted, one way or the other, a more substantial advantage would be seen.  While implicit reverse shows a small performance advantage, this is not a fundamental advantage, but rather an artifact of the implementation.  The reverse mode package allows direct specification of the vector Jacobian product, whereas the forward mode package requires deconstructing and reconstructing the dual numbers (with the Jacobian vector product in between). The reconstruction involves allocation and copying that adds a little extra overhead, particularly for larger problems.

Direct forward and direct reverse have very different performance. As noted, reverse mode is particularly problematic for solvers since taping iterations incurs a large computation and memory penalty.  Direct forward mode, on the other hand, doesn't require storing much extra.  We simply retain the current dual numbers at each operation in the analysis.  In fact, for the very smallest problem, direct forward is actually faster than the implicit methods.  This is not surprising because for a problem with only 2 states, very few solver iterations are required, and so treating the unrolled iterations as explicit operations can potentially be faster.

Finite differencing generally sits between direct forward and direct reverse.  However, as the problem becomes larger the computation time for direct forward grows more rapidly, eventually yielding a slower approach.  This is because we are propagating dual numbers through an increasing number of solver iterations.

The differences are sometimes hard to appreciate on a log scale, and so in \cref{tab:nonjactime} we report the time from the largest problem size, $n = 128$, for each case.  Applying direct reverse to a problem with this many solver iterations is folly, but is done to complete the table, coming it at nearly 45 minutes.  Direct forward AD and central finite differencing were comparable at around 5--7 seconds.   The implicit methods offer over an order of magnitude in speed increase, dropping execution time to one or two tenths of a second.  The adjoint in this case is 50 times faster than finite differencing or 70 times faster than direct forward-mode AD.

\begin{table}[htb]
\centering
\caption{Time in seconds to compute Jacobian for the nonlinear problem with $128$ inputs and $128$ output.}
\label{tab:nonjactime}
\begin{tabular}{@{}ccccc@{}}
\toprule
Direct & Direct & Implicit & Implicit & Central \\
Forward & Reverse & Forward & Reverse & Difference \\
\midrule
6.63 & 2694 & 0.246  & 0.0952 & 4.72 \\
\bottomrule
\end{tabular}
\end{table}

From a developer perspective enabling this speed up is a one-line change from
\\
\jlinl{y = solve(x)}
\\
to
\\
    \jlinl{y = implicit(solve, residual, x)}
\\
where \jlinl{r = residual(x, y)} returns the result of \cref{eq:rosen}, and is a function that is defined anyway for the solve function.  In a normal engineering analysis, explicit calculations would be included before and after the solver.  None of those require any changes as AD automatically propagates the derivatives as usual.  If additional implicit calculations were required we would overload those solvers in the same manner.

From a user perspective there is no change.  One would compute the Jacobian, for example using:
\\
    \jlinl{J = ForwardDiff.jacobian(func, x)}
\\
where \jlinl{func} is some function that contains this nonlinear solve as well as any other operations, for any set of desired design variables or outputs.  The adjoint  occurs automatically without any effort from the user.

\subsection{Unsteady Adjoint for Implicit Differential Equations}

For ODEs solved with an implicit time stepping method we have two mechanisms to speed up the derivative computation.  First, we use implicit differentiation across the solver that occurs at each time step.  Second, we apply reverse-mode AD across each time step separately and analytically propagate derivatives between time steps.  As problems grow in size, this allows us to keep the tape
small and avoid the issues associated with growing memory requirements.

As an example, we consider a 2D heat transfer analysis on a thin plate with convection and radiation.\footnote{Motivated by a similar problem from Mathworks: \url{https://www.mathworks.com/help/pde/ug/nonlinear-heat-transfer-in-a-thin-plate.html}}  The PDE describing the spatial and temporal temperature ($T$) variation from convective and radiative heat transfer is:
\begin{equation}
\frac{\partial T}{\partial t} = \alpha \nabla^2 T - \beta (T - T_a) - \gamma  (T^4 - T_a^4)
\end{equation}
where $\alpha, \beta$, and $\gamma$ are constants, and $T_a$ is the ambient temperature.  Note that we have combined various physical constants (e.g., thermal conductivity of the plate, density of the plate) into the macro-constants $\alpha, \beta$, and $\gamma$ as the physics details are not important for our purposes.

We now discretize the problem with a Cartesian grid, and apply the above governing equation at each node as follows:
\begin{equation}
\begin{aligned}
\frac{d T_{i, j}}{d t} = & \alpha \left(   \frac{T_{i+1, j} - 2 T_{i, j} + T_{i-1, j}}{\Delta x_{i,j}^2} \right.\\
& \quad \left. + \frac{T_{i, j+1} - 2 T_{i, j} + T_{i, j-1}}{\Delta y_{i, j}^2}\right) \\
&\quad - \beta (T_{i, j} - T_a) - \gamma  (T_{i, j}^4 - T_a^4)
\end{aligned}
\end{equation}
With an $n \times n$ grid we have now converted the PDE to $n^2$ ODEs.

We use a square plate, with uniform grid spacing in both directions: $\Delta x_{i, j} = \Delta y_{i, j} = \delta$.  We can then simplify the ODEs as:
\begin{equation}
\begin{aligned}
\frac{d T_{i, j}}{d t} = &\frac{\alpha}{\delta} \left( T_{i+1, j} + T_{i-1, j} + T_{i, j+1} + T_{i, j-1} - 4 T_{i, j}\right) \\
&\quad - \beta (T_{i, j} - T_a) - \gamma (T_{i, j}^4 - T_a^4)
\end{aligned}
\end{equation}

As boundary conditions we specify all edges as insulated, except the bottom edge where we control the temperature.  This means that the left, right, and top edges have a Neumann boundary condition $dT/dx = 0$ or $dT/dy = 0$.  For example, at the left boundary, using our discretization: $T_{i, 1} = T_{i, 2}$ (where we have chosen a matrix layout: rows correspond to different vertical positions on the plate).  Along the bottom row we specify the temperature along all $n$ grid points (a Dirichlet boundary condition), and allow this temperature to change in time as well.  Since the temperature is known on all boundaries, an $n \times n$ grid produces $(n-2) \times (n-2)$ states and the same number of ODEs.  The temperature on the bottom row creates $n \times n_t$ control inputs, where $n_t$ is the number of time steps.

We setup the problem with the following constants: $\alpha = 1.16\times 10^{-4}, \beta = 5.78\times10^{-5}, \gamma = 1.64\times 10^{-12}, T_a = 300$, and specify a temperature gradient along the bottom row from 1000 K on the bottom left boundary to 600 K on the bottom right.  While we have the ability to allow this to change in time, for simplicity we keep the control inputs the same over time. The spacing $\delta$ varies as we change the number of grid points.  We start from the coarsest possible mesh: $3 \times 3$ and increase the mesh density up to $19 \times 19$ where possible.  We simulate across 5,000 seconds, approximately the time it takes to reach steady-state with constant control inputs, and use only 100 time steps with an implicit Euler method.  As output, we return the temperature at the upper left corner of the plate.  Thus, as we vary $n$ from 3 to 19, we have $(n-2)^2$ states, $n \times 1000$ inputs, and $1$ output.  Choices in defining this problem are meant to mimic common engineering optimal control problems where there are many inputs (a set of inputs that can vary at each time step) and a few outputs (typically integral quantities, or final time-step quantities).

\begin{figure}[htbp]
\centering
\includegraphics[width=3.0in]{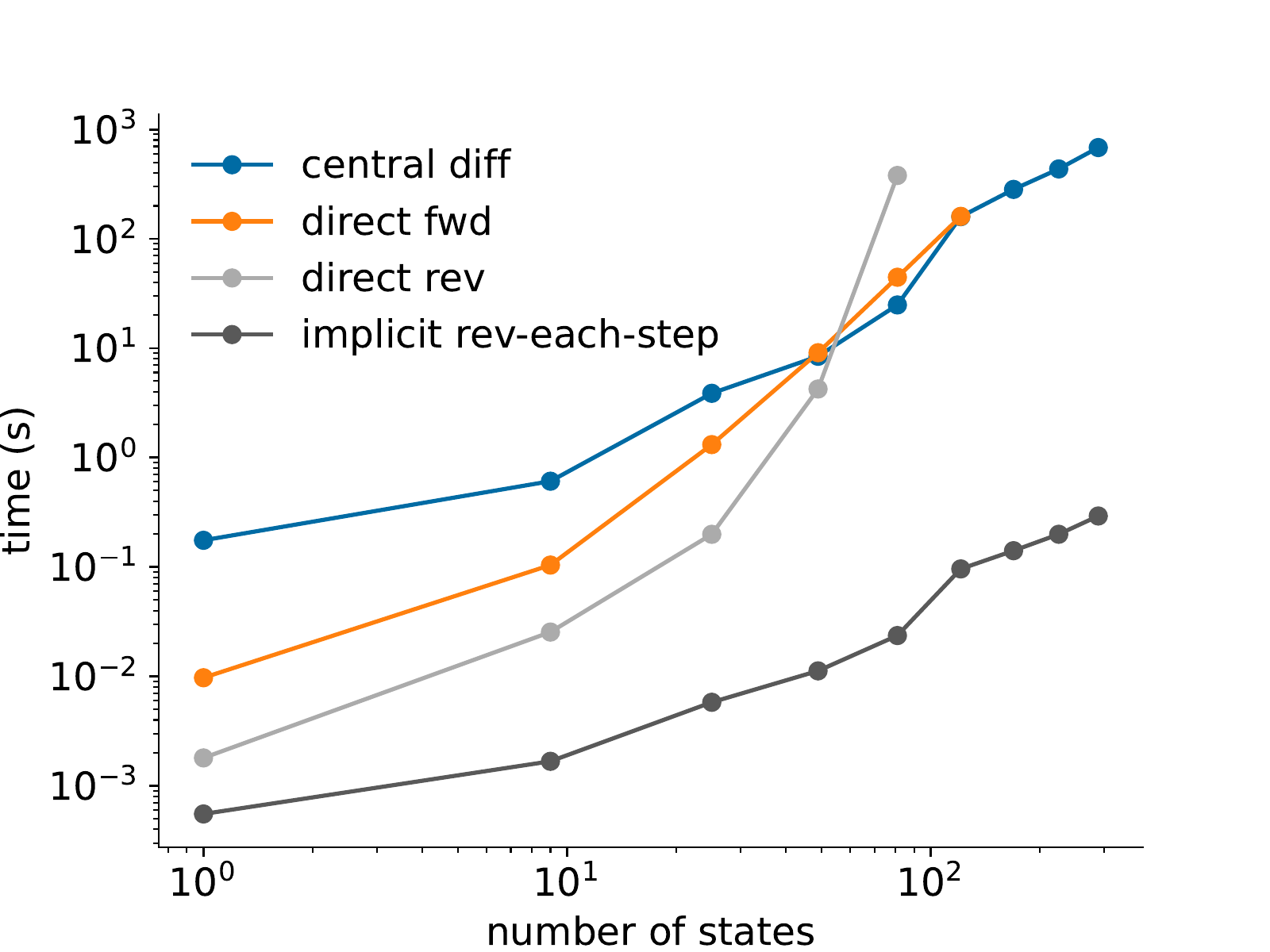}
\caption{Time to compute gradient for the implicit ODE problem as a function of the number of states for finite differencing, forward AD, reverse AD, and reverse implicit AD applied over each time step.}
\label{fig:implicit}
\end{figure}

In \cref{fig:implicit} we plot the time to compute the gradient as a function of the number of states $(n-2)^2$.  In the case of finite-differencing we did not actually finite difference, since the computational cost was quite large, but rather recorded the time for one function call then multiplied by the number of inputs (plus one for evaluating at the current point).  There would be a very small additional overhead in actually finite differencing, thus this is actually a minimum cost for one-sided finite differencing.

For small problems finite differencing is the slowest, with forward AD an order of magnitude faster.  Reverse mode AD is much faster than either, as expected for this problem with many inputs and few outputs.  However, as the number of states increase, the number of iterations required by the solver at each iteration increases also.  Forward mode AD slows down with the extra propagation of dual numbers, but reverse mode AD slows down even more significantly as the computational tape grows rapidly.  Very quickly these methods become intractable.  The implicit reverse mode, however, maintains a much lower computational cost, three orders of magnitude faster than the next fastest method for the largest problem.  For the largest problem run with direct reverse mode, 81 states, the time for each method is reported in \cref{tab:impodetime}.

\begin{table}[htb]
\centering
\caption{Time in seconds to compute gradient for implicit ODE with 81 states.}
\label{tab:impodetime}
\begin{tabular}{@{}cccc@{}}
\toprule
Central Diff & Direct Fwd & Direct Rev & Rev Each Step \\
\midrule
12.4 & 45 & 380 & 0.02 \\
\bottomrule
\end{tabular}
\end{table}

For problems that favor reverse mode, the two modifications shown here offer major improvements.  We avoid the memory issues of a growing computational tape (both from time steps and from internal solver iterations), while also avoiding the unnecessary cost of propagating derivatives through solver iterations.

\subsection{Reverse Mode for Explicit Differential Equations}

For explicit differential equations there are no internal solver iterations to exploit with an adjoint implementation.  But, as discussed,  in \cref{sec:explicitode}, we may still be able to speed up the explicit process by applying reverse-mode AD across each time separately to slow down the growing memory requirements.

We consider the same problem as the previous section, but use the explicit Tsitouras 5/4 Runge-Kutta method \cite{Tsitouras2011}.  Because the problem is explicit we need to increase the number of time steps.  We use 1,000 time steps as that number of time steps allowed the simulation to remain stable across all cases.

\begin{figure}[htbp]
\centering
\includegraphics[width=3.0in]{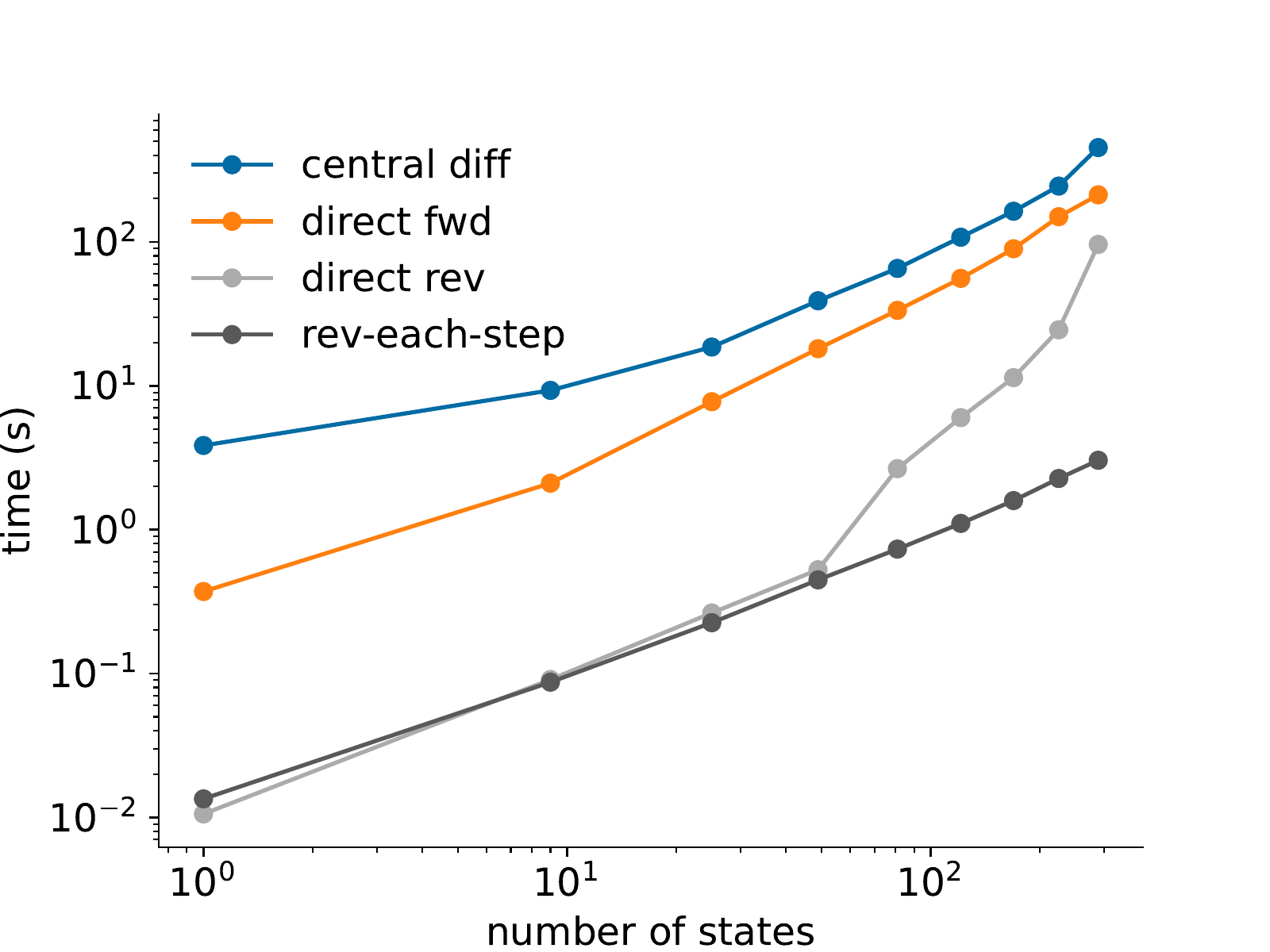}
\caption{Time to compute gradient for the explicit ODE problem as a function of the number of states for finite differencing, forward AD, reverse AD, and reverse AD applied over each time step.}
\label{fig:explicit}
\end{figure}

In \cref{fig:explicit} we plot the time to compute the gradient as a function of the number of states $(n-2)^2$.

As expected, for this problem with many inputs and few outputs, reverse-mode is much faster than forward mode, and finite differencing is intractable.  The relative cost difference doesn't change much with the number of states (except for the two reverse mode methods as the problem becomes large).  This is expected since the number of outputs is fixed, and while the number of inputs is not fixed, it grows much slower (proportional to the square root of the number of states).

For a small number of states, and thus lower complexity in the tape, there is no significant difference between recording a reverse-mode tape across the entire ODE solve versus recording a tape across just one time step.  However, as the number of states grows, and thus the complexity of the tape grows, the one-time-step approach shown here offers significant speed improvements.  Actually, the improvements are even greater than shown as the values only include the time to compute the gradient.  The time to compile the tape is not included since that is one time cost.  However, that one-time cost is \emph{far} greater for the long tape across the entire ODE.
For the last case, the largest number of states, the times are reported in \cref{tab:expodetime}.

\begin{table}[htb]
\centering
\caption{Time in seconds to compute gradient for explicit ODE with 289 states.}
\label{tab:expodetime}
\begin{tabular}{@{}cccc@{}}
\toprule
Central Diff & Direct Fwd & Direct Rev & Rev Each Step \\
\midrule
453 & 212 & 96 & 3 \\
\bottomrule
\end{tabular}
\end{table}

The relative benefits of using this approach are of course problem dependent. But the salient feature that would favor the record-one-time-step approach is tape complexity, which is exacerbated by adding states, increasing the number of time steps, using a more complex ODE method, etc.

\section{Conclusions}

By defining new algorithmic differentiation rules, using implicit differentiation, we can automate the procedure of applying an adjoint.  This approach works for any general nonlinear function, and once implemented, allows any user to apply adjoints with minimal effort.

Many problems have structure, or analytic solutions, that we can leverage, rather than treating all problems as general nonlinear.  Differential equations have a regular structure to which we can apply the automated adjoint formulation in a memory-efficient manner.

Even for explicit differential equations, which have no iterative solvers to leverage, can still be sped up with related techniques.  We record a tape across only one time step, then analytically propagate derivatives in reverse, allowing for efficient derivative computation without forming long tapes.

Across multiple demonstration problems we see that these methods can offer large speed ups, often an order of magnitude or more, with minimal coding changes.  Furthermore, they allow for interfacing with external codes, and other optimization frameworks to offer flexibility while still computing accurate derivatives.

\backmatter

\bmhead{Funding}

This material is based upon work supported by the U.S. Department of Energy, Office of Science, Office of Advanced Scientific Computing Research under Award Number DE-FOA-0002717.
This article was partial developed based upon funding from the Alliance for Sustainable Energy, LLC, Managing and Operating Contractor for the National Renewable Energy Laboratory for the U.S. Department of Energy.

\bmhead{Replication of Results}
A full implementation of the methodology in the Julia programming language, with documentation, is available in an open-source repository called ImplicitAD: \url{https://github.com/byuflowlab/ImplicitAD.jl}.  All of the results will be stored in an ``examples'' subdirectory of this same repository.

\bibliography{impref.bib}

\end{document}